# DISCUSSION OF "BREAKDOWN AND GROUPS" BY P. L. DAVIES AND U. GATHER

By Hannu Oja

*University of Jyväskylä*

**1. Breakdown, equivariance and invariance.** The authors are to be congratulated for their excellent paper, which nicely clarifies the role of equivariance in finding upper bounds for the breakdown points of functionals. The breakdown point approach, with upper bounds showing how far one can go, has achieved great success in the univariate and multivariate location, scale, scatter and regression estimation problems. The authors justifiably argue that this is due to the fact that the acceptable, well-behaved estimates in these contexts have natural equivariance properties. In constructing reasonable estimates and test statistics, one therefore considers statistics satisfying certain conditions (invariance, equivariance, unbiasedness, consistency, etc.). If there are no restrictions, the upper bound is one as the breakdown point (using the common definition) of a "stupid" constant functional, for example, is one.

The paper is clearly written with several illustrative examples. The constructive proof of the main Theorem 3.1 illustrates how one can concretely break down an equivariant estimate:

1. Pick a transformation $g$ corresponding to the set with the supremum probability mass in (3.3).
2. Apply the transformation $g$ or $g^{-1}$ repeatedly to contaminate a (random) half of the data outside the set with the supremum probability mass.

In the one-sample location problem with sample size $n = 2k$, for example, the translation equivariance of a location estimate $T(x_1, \ldots, x_n)$ means that

$$T(x_1 + c, \ldots, x_k + c, x_{k+1}, \ldots, x_n) - T(x_1, \ldots, x_k, x_{k+1} - c, \ldots, x_n - c) = c$$

and, consequently, the estimate can be broken either by repeatedly shifting the first half of the data by $+c$ or by repeatedly shifting the second half of the data by $-c$.

---

The theory thus yields upper bounds for the breakdown points of the affine equivariant univariate location and scale functionals but does not say anything about *affine invariant* skewness and kurtosis statistics, for example. Clearly the invariant classical skewness statistic

$$b_1 = \frac{((1/n)\sum(x_i - \bar{x})^3)^2}{((1/n)\sum(x_i - \bar{x})^2)^3}$$

does not break down with one outlying observation, although affine equivariant second and third central moments both do. For a single outlier going to infinity, the central moments move "beyond all bounds" but $b_1$ converges to a constant $(n-2)^2/(n-1)$. As this limit is not data dependent at all, the contaminated statistic $b_1$ does not convey any information on the original data points. Is this a breakdown? Another strange example is the estimation problem for the parameters of the linear predictor $\theta_0 + \theta'x$ in the generalized linear model. Again, for $n = 2k$, an equivariant estimate of $\theta$ satisfies

$$\hat{\theta}\left(\begin{pmatrix} c \cdot x_1 \\ y_1 \end{pmatrix}, \ldots, \begin{pmatrix} c \cdot x_k \\ y_k \end{pmatrix}, \begin{pmatrix} x_{k+1} \\ y_{k+1} \end{pmatrix}, \ldots, \begin{pmatrix} x_n \\ y_n \end{pmatrix}\right)$$
$$= \frac{1}{c}\hat{\theta}\left(\begin{pmatrix} x_1 \\ y_1 \end{pmatrix}, \ldots, \begin{pmatrix} x_k \\ y_k \end{pmatrix}, \begin{pmatrix} (1/c) \cdot x_{k+1} \\ y_{k+1} \end{pmatrix}, \ldots, \begin{pmatrix} (1/c) \cdot x_n \\ y_n \end{pmatrix}\right)$$

and the estimate can be moved beyond all bounds *or* to zero by repeatedly multiplying half of the data by $c$ or by $1/c$. The estimate then seems to become uninformative. Is something wrong with the definitions of the breakdown and the breakdown point? What do we really mean when we say that a breakdown occurs?

**2. When does the breakdown occur?** Since the early notions by Hampel (1971), the concept of breakdown point has been widely discussed and further developed by several contributors. For considering and comparing different approaches we adopt the following notation. Let $X = (x_1, \ldots, x_n)$ be an original "true" sample of size $n$ lying in the sample space $\mathcal{X}$. The statistic ("estimate") considered is denoted by $T(X)$ with possible values in $\mathcal{T} = \{T(X) : X \in \mathcal{X}\} \subset \mathbb{R}^p$. We say that a point $t$ is interior to $\mathcal{T}$ if it belongs to $\mathcal{T}$ and there is a neighborhood of $t$ which contains only points of $\mathcal{T}$. A point $t \in \mathbb{R}^p$ is exterior to $\mathcal{T}$ if it does not belong to $\mathcal{T}$, and if there exists a neighborhood of $t$ which contains no points of $\mathcal{T}$. Finally, $t$ is called a boundary point of $\mathcal{T}$ if $t$ is neither interior nor exterior to $\mathcal{T}$. Often $\mathcal{T} = \mathbb{R}^p$ and then there are no boundary points.

We next construct a contaminated sample. Let $S = (s_1, \ldots, s_n)$ be a vector of zeros and ones indicating the contamination, and $Y = (y_1, \ldots, y_n)$, also in $\mathcal{X}$, a sample of "outliers." The contaminated sample then consists of the observations $(1 - s_i)x_i + s_i y_i$, $i = 1, \ldots, n$. The number of outlying or alien observations is accordingly $s = \sum s_i$. The contaminated value of the estimate is then

$$T(X, Y, S) = T((1 - s_1)x_1 + s_1 y_1, \ldots, (1 - s_n)x_n + s_n y_n).$$

The breakdown of the estimate is most often defined as follows.

DEFINITION 1. $T$ breaks down at $X$ with $s$ outliers if there exist a sequence $(Y_m)$ in $\mathcal{X}$ and $S$ with $\sum s_i = s$ such that $\|T(X, Y_m, S)\| \to \infty$.

The breakdown point then gives the smallest fraction of outliers $(s/n)$ that suffices to "*drive the estimate beyond all bounds.*" According to this definition, a constant estimate $(T(X) = t_0)$ can never be broken down. Note also that the scale estimate and scatter matrix estimate are usually thought to break down also if they converge to a boundary point of $\mathcal{T}$ (scale estimate converges to zero and the smallest eigenvalue of the scatter matrix estimate converges to zero). In the paper this is taken care of with a suitably chosen pseudometric; see, for example, Section 4.2. Another possibility is to give a new definition:

DEFINITION 2. $T$ breaks down at $X$ with $s$ outliers if there exist a sequence $(Y_m)$ in $\mathcal{X}$ and $S$ with $\sum s_i = s$ such that either (i) $\|T(X, Y_m, S)\| \to \infty$ or (ii) $T(X, Y_m, S) \to t_0$ where $t_0$ is a boundary point of $\mathcal{T}$.

If the constant functional $T(X) = t_0$ does not depend on $X$, $\mathcal{T} = \{t_0\}$ and $t_0$ is also a boundary point. Therefore $T$ breaks down for all $S$. Note that the boundary point $t_0$ in the definition may depend on $X$, however. In the simple regression example suggested by the referee and analyzed in Section 6, the statistic $T(P_n)$ has values in $\mathcal{T} = [-n, n]$ and it breaks down (in the sense of Definition 2) if $\sum s_i = 1$.

Genton and Lucas (2003) take a different viewpoint and argue that a crucial property of an estimator $T(X, Y, S)$ is that it takes different values for different values of $X \in \mathcal{X}$ and that the breakdown occurs if this property is lost. In this spirit one can say that:

DEFINITION 3. $T$ breaks down with $s$ outliers if there exist a sequence $(Y_m)$ in $\mathcal{X}$ and $S$ with $\sum s_i = s$ such that either (i) $\|T(X, Y_m, S)\| \to \infty$, for all $X \in \mathcal{X}$, or (ii) $T(X, Y_m, S) \to t_0 \in \mathbb{R}^p$, for all $X \in \mathcal{X}$.

In this definition, it is remarkable that the interior or boundary point $t_0$ is not allowed to depend on $X$ at all. This definition solves the problem with the classical skewness statistic; $b_1$ can be made to break down with a single extreme outlier. I wonder whether the techniques and results in the paper by Davies and Gather could be expanded to cover this definition also.

Genton and Lucas (2003) seem in fact to be still more permissive and say that $T(X, Y, S)$ breaks down if

$$\mathcal{T} \cap \left\{ \lim_m T(X, Y_m, S) : X \in \mathcal{X} \right\}$$

collapses to a finite set; an empty set and a singleton $\{t_0\}$ are then special cases. Given a continuum of values of $X$, one expects a continuum of possible values of the estimate. In the linear predictor estimation problem this definition implies that the breakdown point of an equivariant estimate of $\theta$ is at most one half.

All the approaches described above work with the worst possible scenario represented by a strategically chosen sequence of the sets of outlying observations $(Y_m)$. In practice, the observed contaminated value of the estimate $T(X, Y, S)$ is in $\mathcal{T}$, however, and not a boundary point, and one can ask whether the estimate still conveys useful information about the true data cloud or not. Then, instead of speculating about the sequences $(Y_m)$, one may consider the set of possible values of $T(X, Y, S)$ for all choices of $Y \in \mathcal{X}$. With (at most) $s$ outliers, the set of possible values of $T(X, Y, S)$ is

$$\mathcal{T}_s(X) := \left\{ T(X, Y, S) : Y \in \mathcal{X}, \sum s_i = s \right\}.$$

Then clearly

$$\{T(X)\} = \mathcal{T}_0(X) \subset \mathcal{T}_1(X) \subset \mathcal{T}_2(X) \subset \cdots \subset \mathcal{T}$$

and the value of the estimate is totally determined by $s$ outliers if $\mathcal{T}_s(X) = \mathcal{T}$. More generally, we can define that:

DEFINITION 4. $T$ breaks down with $s$ outliers if the set $\mathcal{T}_s(X)$ does not depend on $X$.

Note that if $T$ is affine equivariant/invariant, then also $\mathcal{T}_s(X)$ is affine equivariant/invariant. Assume next that the observed value of $T(X, Y, S)$ is $t$. If we knew the maximum number of outliers in the data set, say $s$, but $S$ and $Y$ are unknown, the observed event $\mathcal{T}_s(X) \ni t$ may still be informative. In the univariate location case with $n = 2k - 1$ and $\mathcal{X} = \mathbb{R}^n$, $\mathcal{T}_s(X) = \mathbb{R}$ for the sample mean if $s > 0$. But for $\mathcal{X} = [0, \infty)^n$, for example, the breakdown point of the mean is one as $\mathcal{T}_{s-1}(X) \ni t \iff x_{(1)} \leq n \cdot t$. For the sample median, the event

$$\mathcal{T}_s(X) \ni t \iff x_{(k-s)} \leq t \leq x_{(k+s)}, \qquad s = 0, \ldots, k,$$

is clearly data dependent and therefore carries information about the data cloud.

DEPARTMENT OF MATHEMATICS
 AND STATISTICS
UNIVERSITY OF JYVÄSKYLÄ
P.O. BOX 35
FIN 40351 JYVÄSKYLÄ
FINLAND
E-MAIL: ojahannu@maths.jyu.fi
URL: www.maths.jyu.fi/~ojahannu